\documentclass[letter,10pt]{amsart}
\usepackage[english]{babel}
\usepackage{multicol}
\usepackage{graphicx}
\usepackage{amssymb}
\usepackage[all]{xy}
\usepackage{amsmath,amsfonts}
\usepackage{enumerate}
\usepackage{array}
\usepackage{tikz}
\input xy
\xyoption{all}

\newtheorem{teo}{Theorem}[section]
\newtheorem{cor}[teo]{Corollary}
\newtheorem{lema}[teo]{Lemma}
\newtheorem{prop}[teo]{Proposition}
\newtheorem{defi}[teo]{Definition}

\newtheorem{exam}[teo]{Example}

\newtheorem{Rem}[teo]{Remark}

\numberwithin{equation}{section}

\title{The third homology of $SL_{2}$ of real quadratically closed fields.}
\author{Rodrigo Cuitun Coronado}
\address{School of Mathematics and Statistics, University College Dublin}
\email{rodrigo.cuituncoronado@ucdconnect.ie}
\date{\today}

\begin{document}


\begin{abstract}
For a real closed field $\mathbf{R}$, we use the theory of the refined Bloch group to give a new short proof of the isomorphisms $H_{3}(SL_{2}(\mathbf{R}),\mathbb{Z})\cong K_{3}^{\mathrm{ind}}(\mathbf{R})$ and $H_{3}(SL_{n}(\mathbf{R}),\mathbb{Z})\cong H_{3}(SL_{2}(\mathbf{R}),\mathbb{Z})\oplus K_{3}^{M}(\mathbf{R})^{0}$ for $n\geq3$.
\end{abstract}

\maketitle


\section{Introduction.}

In the paper \emph{Homology of classical Lie groups made discrete III} (\cite{ChinHanSah}) the author, Chih-Han Sah, proves the injectivity of the map $H_{3}(SL_{2}(F),\mathbb{Z})\rightarrow H_{3}(SL_{3}(F),\mathbb{Z})$ and $H_{3}(SL_{3}(F),\mathbb{Z})\cong H_{3}(SL_{2}(F),\mathbb{Z})\oplus K_{3}^{M}(F)^{0}$ for $F=\mathbb{R}$, $\mathbb{H}$ or $F$ algebraically closed (\cite[Theorem 3.0]{ChinHanSah} and \cite[Proposition 2.15]{ChinHanSah} respectively). He states, without giving details, that the result is also valid for any real closed field \cite[Section 3]{ChinHanSah}.

In the present article we use the theory of the refined Bloch group (\cite{ArticleBlochWignerComplex}) to give a new proof of this result for real closed fields $\mathbf{R}$, and we deduce the structure of $H_{3}(SL_{n}(\mathbf{R}),\mathbb{Z})$, $n\geq3$, in terms of $K$-theory. In fact, all the results are valid for the larger class of real quadratically closed fields (See section $2$ below).


\subsection{Some notation}

For a field $F$, we let $F^{\times}$ denote the group of units of $F$. For $x\in F^{\times}$ we will let $\langle x\rangle\in F^{\times}/(F^{\times})^{2}$ denote the corresponding square class. Let $R_{F}$ denote the integral group ring $\mathbb{Z}\left[F^{\times}/(F^{\times})^{2}\right]$ of the group $F^{\times}/(F^{\times})^{2}$. We will use the notation $\langle\langle x\rangle\rangle$ for the basis element, $\langle x\rangle-1$, of the augmentation ideal $\mathcal{I}_{F}$ of $R_{F}$.


\section{Real quadratically closed fields.}

Recall that an \emph{ordered field} $(\mathbf{R},\geq)$ is \emph{real closed} if it satisfies the following two properties:

\begin{itemize}
\item[(1)]Any positive element $x>0$ in $\mathbf{R}$ has a square root in $\mathbf{R}$; i.e. there exist $a\in \mathbf{R}\setminus\{0\}$ such that $x=a^{2}$.
\item[(2)] Any odd degree polynomial with coefficients in $\mathbf{R}$ has a root in $\mathbf{R}$.
\end{itemize}

In fact, our main results will apply to any ordered field satisfying property $(1)$. For convenience, in this paper we will refer to such fields as \emph{real quadratically closed fields}

\begin{exam}\label{Fieldconstructiblerealnumbers}
The field of (straightedge-and-compass) constructible real numbers is real quadratically closed but not real closed.
\end{exam}

We will frequently use the following easily verified fact

\begin{lema}
Let $(\mathbf{R},\geq)$ be an ordered field and $\mathbf{R}_{>0}$ denote the  group of positive elements in $\mathbf{R}$ under the multiplication of $\mathbf{R}$.

\begin{itemize}
\item[(i)]$\mu_{\mathbf{R}}=\mu_{2}$ and $\mathbf{R}^{\times}=\mu_{2}\oplus\mathbf{R}_{>0}$.
\item[(ii)] If $\mathbf{R}$ is real quadratically closed then $\mathbf{R}_{>0}=(\mathbf{R}_{>0})^{2}$
\end{itemize}
\end{lema}

\begin{flushright}
    $\Box$
\end{flushright}

Since $\mathbf{R}_{>0}=(\mathbf{R}_{>0})^{2}$, it follows that $\mathbf{R}^{\times}/(\mathbf{R}^{\times})^{2}=\mu_{2}$ and hence $R_{\mathbf{R}}=\mathbb{Z}[\mu_{2}]=\mathbb{Z}[\langle-1\rangle]$. Futhermore we will need the following observation:

\begin{cor}\label{corollaryrealquadraticallyclosedfields}If $\mathbf{R}$ is real quadratically closed then $\bigwedge^{2}_{\mathbb{Z}}(\mathbf{R}^{\times})\cong\bigwedge^{2}_{\mathbb{Z}}(\mathbf{R}_{>0})$ which is $2$-divisible
\end{cor}

\begin{flushright}
    $\Box$
\end{flushright}


\section{Review of Bloch groups of fields.}

In this section we will recall the definition and applications of the classical pre-Bloch group $\mathcal{P}(F)$ and the refined pre-Bloch group $\mathcal{RP}(F)$.


\subsection{Classical Bloch Group $\mathcal{B}(F)$.}

For a field $F$, with at least $4$ elements, the \emph{pre-Bloch group} or \emph{Scissors congruence group}, $\mathcal{P}(F)$, is the group generated by the elements $[x]$, with $F^{\times}\setminus\{1\}$, subject to the relations

\begin{center}
$R_{x,y}:\; 0=[x]-[y]+[y/x]-[(1-x^{-1})/(1-y^{-1})]+[(1-x)/(1-y)]$, \;\; $x\not=y$.
\end{center}

Let $S^{2}_{\mathbb{Z}}(F^{\times})$ denote the group

\begin{center}
$\frac{F^{\times}\otimes_{\mathbb{Z}}F^{\times} }{\langle x\otimes y +y\otimes x\; |\;x,y\in F^{\times}\rangle}$
\end{center}

and denote by $x\circ y$ the image of $x\otimes y$ in $S^{2}_{\mathbb{Z}}(F^{\times})$. The map

\begin{center}
$\lambda: \mathcal{P}(F)\rightarrow S^{2}_{\mathbb{Z}}(F^{\times})$,  \;\; $[x]\mapsto (1-x)\circ x$
\end{center}


is well-defined, and the \emph{Bloch group} of $F$, $\mathcal{B}(F)\subset\mathcal{P}(F)$, is defined to be the kernel of $\lambda$.

The Bloch group $\mathcal{B}(F)$ of a general field $F$ is of interest because of the following result of Suslin on $K_{3}^{\mathrm{ind}}$:

\begin{teo}\cite[Theorem 5.2]{Suslin1}\label{Suslinmaintheorem}
Let $F$ be an infinite field, then there is a short exact sequence

\begin{center}
$\xymatrix{0\ar[r]&\mathrm{Tor}_{1}^{\mathbb{Z}}(\widetilde{\mu_{F},\mu_{F}})\ar[r]&K_{3}^{\mathrm{ind}}(F)\ar[r]&\mathcal{B}(F)\ar[r]&0}$
\end{center}

where $\mathrm{Tor}_{1}^{\mathbb{Z}}(\widetilde{\mu_{F},\mu_{F}})$  is the unique nontrivial extension of $\mathrm{Tor}_{1}^{\mathbb{Z}}(\mu_{F},\mu_{F})$ by $\mathbb{Z}/2$ when $\mathrm{Char}(F)\not=2$ (and $\mathrm{Tor}_{1}^{\mathbb{Z}}(\widetilde{\mu_{F},\mu_{F}})=\mathrm{Tor}_{1}^{\mathbb{Z}}(\mu_{F},\mu_{F})$ if $\mathrm{Char}(F)=2$) .
\end{teo}

\begin{flushright}
    $\Box$
\end{flushright}


\subsection{The refined Bloch Group $\mathcal{RB}(F)$.}

The \emph{refined pre-Bloch group} $\mathcal{RP}(F)$, of a field $F$ which has at least $4$ elements, is the $R_{F}$-module with generators $[x]$, $x\in F^{\times}$ subject to the relations $[1]=0$ and

\begin{center}
$S_{x,y}:\; 0=[x]-[y]+\langle x\rangle[y/x]-\langle x^{-1}-1\rangle[(1-x^{-1})/(1-y^{-1})]+\langle1-x\rangle[(1-x)/(1-y)]$, \;\; $x,y\not=1$.
\end{center}

From the definitions, it follows that $\mathcal{P}(F)=(\mathcal{RP}(F))_{F^{\times}}$. Let $\Lambda=(\lambda_{1},\lambda_{2})$ be the $R_{F}$-module homomorphism

\begin{center}
$\mathcal{RP}(F)\rightarrow \mathcal{I}_{F}^{2}\oplus S^{2}_{\mathbb{Z}}(F^{\times})$
\end{center}

where $\lambda_{1}:\mathcal{RP}(F)\rightarrow \mathcal{I}^{2}_{F}$ is the map $[x]\mapsto \langle\langle1-x\rangle\rangle\langle\langle x\rangle\rangle$ and $\lambda_{2}$ is the composite

\begin{center}
$\xymatrix{ \mathcal{RP}(F)\ar[r]&\mathcal{P}(F)\ar[r]^{\lambda} &S^{2}_{\mathbb{Z}}(F^{\times})}$,
\end{center}

where $S^{2}_{\mathbb{Z}}(F^{\times})$ has the trivial $R_{F}$-module structure.


The \emph{refined Bloch group} of $F$ is the module

\begin{center}
$\mathcal{RB}(F):=\mathrm{Ker}(\Lambda: \mathcal{RP}(F)\rightarrow \mathcal{I}_{F}^{2}\oplus S^{2}_{\mathbb{Z}}(F^{\times}  )$.
\end{center}

Furthermore, the \emph{refined scissor congruence group} of $F$ is the $R_{F}$-module

\begin{center}
$\mathcal{RP}_{1}(F):=\mathrm{Ker}(\lambda_{1}:\mathcal{RP}(F)\rightarrow\mathcal{I}^{2}_{F})$.
\end{center}

Thus $\mathcal{RB}(F)=\mathrm{Ker}(\lambda_{2}:\mathcal{RP}_{1}(F)\rightarrow S^{2}_{\mathbb{Z}}(F^{\times}))$.

Now for any infinite field $F$, Suslin shows that there is a natural homomorphism $H_{3}(SL(F),\mathbb{Z})\rightarrow K_{3}^{\mathrm{ind}}(F)$, where $SL(F)=\displaystyle\bigcup_{n\geq1}SL_{n}(F)$ and $K_{3}^{\mathrm{ind}}(F)$ is the indecomposable $K_{3}$ of $F$. It follows that there are induced maps $H_{3}(SL_{n}(F),\mathbb{Z})\rightarrow K_{3}^{\mathrm{ind}}(F)$ for all $n$. Hutchinson and Tao show that $H_{3}(SL_{2}(F),\mathbb{Z})\rightarrow K_{3}^{\mathrm{ind}}(F)$ is surjective for any infinite field \cite[Lemma 5.1]{HutchinsonLiqunTao}.


The refined Bloch group is of interest because of the following result on the third homology of $SL_{2}$ over a general field $F$:

\begin{teo}\cite[Theorem 4.3]{ArticleBlochWignerComplex}\label{TheoremBlochWignerComplex}
There is a complex of $R_{F}$-modules

\begin{center}
$\xymatrix{0\ar[r]&\mathrm{Tor}_{1}^{\mathbb{Z}}(\mu_{F},\mu_{F})\ar[r]&H_{3}(SL_{2}(F),\mathbb{Z})\ar[r]& \mathcal{RB}(F)\ar[r]&0}$,
\end{center}

which is exact except at the middle term where the homology is annihilated by $4$.

\end{teo}

\begin{flushright}
    $\Box$
\end{flushright}

Furthermore we have that the following diagram is commutative:

\begin{center}
$\xymatrix{0\ar[r]&\mathrm{Tor}_{1}^{\mathbb{Z}}(\mu_{F},\mu_{F})\ar[r]\ar[d]&H_{3}(SL_{2}(F),\mathbb{Z})\ar[r]\ar[d]& \mathcal{RB}(F)\ar[r]\ar[d]&0\\
0\ar[r]&\mathrm{Tor}_{1}^{\mathbb{Z}}(\widetilde{\mu_{F},\mu_{F}})\ar[r]&K_{3}^{\mathrm{ind}}(F)\ar[r]& \mathcal{B}(F)\ar[r]&0}$,
\end{center}

\begin{Rem}
We will examine the sequence of Theorem \ref{TheoremBlochWignerComplex} below in more detail in the case of real quadratically closed fields.
\end{Rem}


\subsection{The elements $\psi_{i}(x)$ in $\mathcal{RP}(F)$}

For $x\in F^{\times}$, we define the following elements of $\mathcal{RP}(F)$

\begin{center}
$\psi_{1}(x):=[x]+\langle-1\rangle[x^{-1}]$ \;\; and\;\; $\psi_{2}(x):=\left\{
                                                                       \begin{array}{ll}
                                                                         \langle x^{-1}-1\rangle[x]+\langle1-x\rangle[x^{-1}], & \hbox{$x\not=1$;} \\
                                                                         0, & \hbox{$x=1$.}
                                                                       \end{array}
                                                                     \right.
$.
\end{center}

From the definitions of the the elements $\psi_{i}(x)$, we get that $\langle-1\rangle\psi_{i}(-1)=\psi_{i}(-1)$ for $i\in\{1,2\}$. We review some of the  algebraic properties of the elements $\psi_{i}(x)$:

\begin{prop}\cite[Lemma 3.1, Proposition 3.2]{ArticleBloch}\label{algebraicpropertiespsi} Let $F$ be a field. For $i\in\{1,2\}$, we have:

\begin{itemize}
\item[(1)] $\psi_{i}(xy)=\langle x\rangle\psi_{i}(y)+\psi_{i}(x)$, for all $x,y$.
\item[(2)] $\langle\langle x\rangle\rangle\psi_{i}(y)=\langle\langle y\rangle\rangle\psi_{i}(x)$, for all $x,y$.
\item[(3)] $2\psi_{i}(-1)=0$.
\item[(4)] $\psi_{i}(x^{2})=\langle\langle x\rangle\rangle\psi_{i}(-1)$, for all $x$.
\item[(5)] $2\psi_{i}(x^{2})=0$ for all $x$ and if $-1\in (A^{\times})^{2}$, then $\psi_{i}(x^{2})=0$ for all $x$.
\end{itemize}

\end{prop}

\begin{flushright}
    $\Box$
\end{flushright}


\subsection{The constant $C_{F}$ .}

In the section $3.2$ of \cite{ArticleBloch}, it is shown that the elements

\begin{center}
$C(x)=[x]+\langle-1\rangle[1-x]+ \langle\langle1-x\rangle\rangle\psi_{1}(x) \in \mathcal{RP}(F)$
\end{center}

are constant for a field with at least 4 elements (Lemma 3.5); i.e. $C(x)=C(y)$ for all $x,y\in F^{\times}$. Therefore we have the following definition.

\begin{defi}
Let $F$ be a field with at least $4$ elements. We will denote by $C_{F}$ the common value of the expression $C(x)$ for $x\in F\setminus\{0,1\}$; i.e.

\begin{center}
$C_{F}:=[x]+\langle-1\rangle[1-x]+ \langle\langle1-x\rangle\rangle\psi_{1}(x)$ \;\;in\;\; $\mathcal{RP}(F)$.
\end{center}

\end{defi}

One can deduce some properties of the element $C_{F}$

\begin{prop}\label{actiontrivialonelementCF}\cite[Corollary 3.7]{ArticleBloch}
Let $F$ be a field. Then

\begin{itemize}
\item[(1)]$\langle-1\rangle C_{F}=C_{F}$.
\item[(2)]For any field $F$, $C_{F}\in \mathcal{RB}(F)$; i.e. $\Lambda(C_{F})=0$.
\end{itemize}
\end{prop}

\begin{flushright}
    $\Box$
\end{flushright}

Observe that the image of $C_{F}$ in $\mathcal{P}(F)$ is $[x]+[1-x]$ for $x\not=0,1$. We also denote this element of $\mathcal{B}(F)\subset \mathcal{P}(F)$ by $C_{F}$.


\section{The Refined Bloch Group of a real quadratically closed field}

In this section we show that $\langle-1\rangle$ acts trivially on $\mathcal{RP}(\mathbf{R})$ where $\mathbf{R}$ is a real quadratically closed field; i.e. $R_{\mathbf{R}}$ acts trivially on $\mathcal{RP}(\mathbf{R})$.

\begin{lema}\label{psiequalzero}
Let $x\in\mathbf{R} $, if $x>0$ then $\psi_{i}(x)=0$ for $i\in\{1,2\}$ in $\mathcal{RP}(\mathbf{R})$.
\end{lema}

$\textbf{\emph{Proof}}.$ If $x>0$, then $x=a^{2}$ for some $a\in \mathbf{R}$ and hence $\psi_{i}(x)=\langle\langle a\rangle\rangle\psi_{i}(-1)$ by Proposition \ref{algebraicpropertiespsi} (4). Replacing $a$ by $-a$ if necessary, we can suppose $a>0$. Thus $a$ is a square so that $\langle\langle a\rangle\rangle=0$ and hence $\psi_{i}(x)=0$.

\begin{flushright}
    $\Box$
\end{flushright}


\begin{prop}\label{trivialaction}
Let $x\in\mathbf{R} $, if $x>0$ then the action of $\langle-1\rangle$ on $[x]$ is trivial i.e. $\langle-1\rangle[x]=[x]$.
\end{prop}

$\textbf{\emph{Proof}}.$ Since $x>0$, we have that $x=a^{2}$ for some $a\in \mathbf{R}\setminus\{0\}$. Therefore

\begin{center}
$\langle x^{-1}-1\rangle=\langle1-x\rangle\langle x\rangle=\langle1-x\rangle\langle a^{2}\rangle=\langle 1-x\rangle$.
\end{center}

By Lemma \ref{psiequalzero}, we have that

\begin{center}
$0=\psi_{2}(x)=\langle x^{-1}-1\rangle[x]+\langle1-x\rangle[x^{-1}]=\langle1-x\rangle([x]+[x^{-1}])$.
\end{center}

Multiplying the equation $\langle1-x\rangle([x]+[x^{-1}])=0$ by $\langle1-x\rangle$, we get that $[x]+[x^{-1}]=0$, whence we get that $[x]=-[x^{-1}]$. Similarly applying Lemma \ref{psiequalzero} to $\psi_{1}(x)$, we get that $[x]=-\langle-1\rangle[x^{-1}]$. By combining these last two equalities, we get our desired result.

\begin{flushright}
    $\Box$
\end{flushright}

\begin{prop}\label{ElementCAonnegativeelements}
If $x<0$ then $C_{\mathbf{R}}=[x]+[1-x]$ in $\mathcal{RP}(\mathbf{R})$.
\end{prop}

$\textbf{\emph{Proof}}.$ If $x<0$ then $1-x>0$. Hence $\langle\langle1-x\rangle\rangle\psi_{1}(x)=0$ and by Proposition \ref{trivialaction}, we get that $\langle-1\rangle[1-x]=[1-x]$. It follows that $C_{\mathbf{R}}=[x]+[1-x]$ in $\mathcal{RP}(\mathbf{R})$.

\begin{flushright}
    $\Box$
\end{flushright}

\begin{prop}\label{trivialactiononnegatives}
Let $x\in\mathbf{R} $, if $x<0$ then the action of $\langle-1\rangle$ on $[x]$ is trivial i.e. $\langle-1\rangle[x]=[x]$.
\end{prop}

$\textbf{\emph{Proof}}.$ Since $1-x>0$, by Proposition \ref{trivialaction}, we get that $\langle-1\rangle[1-x]=[1-x]$. Therefore from this latter result, Proposition \ref{ElementCAonnegativeelements} and Proposition \ref{actiontrivialonelementCF} $(1)$, we obtain that:

\begin{center}
$0=C_{\mathbf{R}}-\langle-1\rangle C_{\mathbf{R}}=[x]+ [1-x]-\langle-1\rangle[x]-\langle-1\rangle[1-x]=[x]-\langle-1\rangle[x]$.
\end{center}

It follows that $\langle-1\rangle[x]=[x]$.

\begin{flushright}
    $\Box$
\end{flushright}


\begin{cor} \label{RP(R)=P(R)}For any real quadratically closed field $\mathbf{R}$. The natural map $\mathcal{RP}(\mathbf{R})\rightarrow\mathcal{P}(\mathbf{R})$ is an isomorphism.
\end{cor}

$\textbf{\emph{Proof}}.$ $\mathcal{RP}(\mathbf{R})=\mathcal{RP}(\mathbf{R})_{\mathbf{R}^{\times}}\cong \mathcal{P}(\mathbf{R})$ by Proposition \ref{trivialaction} and Proposition \ref{trivialactiononnegatives}

\begin{flushright}
    $\Box$
\end{flushright}

\begin{cor}\label{RP1(R)=P(R)}
For any real quadratically closed field $\mathbf{R}$, the natural map $\mathcal{RP}_{1}(\mathbf{R})\rightarrow\mathcal{P}(\mathbf{R})$ is an isomorphism.
\end{cor}

$\textbf{\emph{Proof}}.$ Note that the map $\lambda_{1}:\mathcal{RP}(\mathbf{R})\rightarrow \mathcal{I}_{\mathbf{R}}^{2}$ given by $[x]\mapsto \langle\langle1-x\rangle\rangle\langle\langle x\rangle\rangle$ is the zero map for all $x\in \mathbf{R}^{\times}$ since either $x>0$ or $1-x>0$. Hence $\mathcal{RP}_{1}(\mathbf{R})=\mathcal{RP}(\mathbf{R})$.

\begin{flushright}
    $\Box$
\end{flushright}


\begin{cor}\label{RB(R)=B(R)}
For any real quadratically closed field $\mathbf{R}$, the natural map $\mathcal{RB}(\mathbf{R})\rightarrow\mathcal{B}(\mathbf{R})$ is an isomorphism.
\end{cor}

$\textbf{\emph{Proof}}.$ $\mathcal{RB}(\mathbf{R})=\mathrm{Ker}(\lambda_{2}:\mathcal{RP}_{1}(\mathbf{R})\rightarrow S^{2}_{\mathbb{Z}}(\mathbf{R}^{\times}))=\mathrm{Ker}(\lambda:\mathcal{P}(\mathbf{R})\rightarrow S^{2}_{\mathbb{Z}}(\mathbf{R}^{\times}) )=\mathcal{B}(\mathbf{R})$ by Corollary \ref{RP1(R)=P(R)}.

\begin{flushright}
    $\Box$
\end{flushright}


\section{The third homology of $SL_{2}(\mathbf{R})$}

We consider the spectral sequence used in section $4$ of \cite{ArticleBlochWignerComplex} to relate $H_{3}(SL_{2}(F),\mathbb{Z})$ to $\mathcal{RB}(F)$.

First we recall two standard subgroups of $SL_{2}(F)$.

\begin{center}
$T:=\left\{\left(
            \begin{array}{cc}
              a & 0 \\
              0 & a^{-1} \\
            \end{array}
          \right)
\;|a\in F^{\times}\right\}$ \;\;\; $B:=\left\{\left(
            \begin{array}{cc}
              a & b \\
              0 & a^{-1} \\
            \end{array}
          \right)
\;|a\in F^{\times}\;, b\in F\right\}$.
\end{center}

We will now examine this spectral sequence in the case of real quadratically closed fields. We will need the following lemma:

\begin{lema}\cite[Lemma 4.1]{ArticleBlochWignerComplex}\label{IsomorphismHomologiesofTandB}If $F$ is an infinite field, then the inclusion $T\rightarrow B$ induces homology isomorphisms

\begin{center}
$H_{k}(T,\mathbb{Z})\cong H_{k}(B,\mathbb{Z})$
\end{center}

for all $k\geq0$.

\end{lema}

\begin{flushright}
    $\Box$
\end{flushright}


\subsection{The Spectral Sequence.}
Let $F$ be a field, and $G=SL_{2}(F)$ act (on the left) on $\mathbb{P}^{1}(F)$ by fractional linear transformations. Let $X_{n}$ be the set of ordered $(n+1)$-tuples $(x_{0},\ldots,x_{n})$ of distinct points of $\mathbb{P}^{1}(F)$. $X_{n}$ is a left $G$-set by diagonal action. Let $L_{n}=\mathbb{Z}[X_{n}]$ and let $d_{n}:L_{n}\rightarrow L_{n-1} $ given by

\begin{center}
$d_{n}(x_{0},\ldots,x_{n})=\displaystyle\sum_{i=0}^{n}(-1)^{i}(x_{0},\ldots,\widehat{x}_{i},\ldots,x_{n})$.
\end{center}

The natural augmentation $\varepsilon:L_{0}\rightarrow \mathbb{Z}$ gives a weak equivalence  $L_{\bullet}\rightarrow \mathbb{Z}$ where $\mathbb{Z}$ is a complex concentrated in degree $0$. It follows that there is a spectral sequence

\begin{center}
$E^{1}_{p,q}=H_{p}(G,L_{q})\Rightarrow H_{p+q}(G,\mathbb{Z})$.
\end{center}

In fact $E^{1}_{0,q}=H_{0}(G,L_{q})=(L_{q})_{G}$.

\begin{prop}\label{TransitivityofSL2} Let $F$ be a field

\begin{itemize}
\item[(1)] $SL_{2}(F)$ acts transitively on $X_{1}$. The stabilizer of $(\infty)$ is the subgroup $B$.
\item[(2)] $SL_{2}(F)$ acts transitively on $X_{2}$. The stabilizer of $(\infty,0)$ is the subgroup $T$.
\end{itemize}

\end{prop}

\begin{flushright}
    $\Box$
\end{flushright}

\begin{cor}\label{Firsttworowsspectralsequence}Let $F$ be an infinite field, then

\begin{itemize}
\item[(1)] $E^{1}_{p,0}\cong H_{p}(G,\mathbb{Z}[G/B])\cong H_{p}(B,\mathbb{Z})\cong H_{p}(T,\mathbb{Z})$.
\item[(2)] $E^{1}_{p,0}\cong H_{p}(G,\mathbb{Z}[G/T])\cong H_{p}(T,\mathbb{Z})$.
\end{itemize}

\end{cor}
\begin{flushright}
    $\Box$
\end{flushright}

\begin{prop}\cite[Section 2.3]{ArticleBlochWignerComplex}\label{IsomorphismXnandZn} Let $F$ be an infinite field. Let

\begin{center}
$\displaystyle\phi(x,y,z):=\left\{
                \begin{array}{ll}
                  (z-x)(x-y)(z-y)^{-1}, & \hbox{$x,y,z\not=\infty$;} \\
                  (y-z)^{-1}, & \hbox{$x=\infty$;} \\
                  z-x, & \hbox{$y=\infty$;} \\
                  x-y, & \hbox{$z=\infty$.}
                \end{array}
              \right.
$.
\end{center}

 Let $Z_{0}=\emptyset$ and for $n\geq1$, let $Z_{n}$ denote the set of ordered $n$-tuples $[z_{1},\ldots,z_{n}]$ of distinct points of $F^{\times}\setminus\{1\}$. Then for $n\geq2$, there is a isomorphism of $R_{F}-$modules

\begin{center}
$E^{1}_{0,n+1}=(L_{n+1})_{G}\longleftrightarrow R_{F}[Z_{n-2}]$
\end{center}

\begin{center}
$(x_{0},\ldots,x_{n})\longmapsto\langle\phi(x_{0},x_{1},x_{2})\rangle\left[\frac{\phi(x_{0},x_{1},x_{3})}{\phi(x_{0},x_{1},x_{2})},\cdots,\frac{\phi(x_{0},x_{1},x_{n})}{\phi(x_{0},x_{1},x_{2})} \right]$
\end{center}

\begin{center}
$(0,\infty,1,z_{1},\ldots,z_{n-2})\longleftarrow[z_{1},\ldots,z_{n-2}]$
\end{center}

\end{prop}

\begin{flushright}
    $\Box$
\end{flushright}


\begin{lema}\cite[Section 4.4]{ArticleBlochWignerComplex}\label{Firstpagespectralsequence}Let $F$ be an infinite field, the $E^{1}$-page of the spectral sequence $E^{1}_{p,q}=H_{n}(G,L_{\bullet})\Rightarrow H_{n}(G,\mathbb{Z})$ has the form

\begin{center}
$\scalebox{0.85}{\xymatrix{&\ar@{-}[ddddddd]&\vdots &\vdots&\vdots&\vdots&\\
&& R_{F}[Z_{2}]\ar[d]^{d^{1}}&R_{F}[Z_{2}]\otimes\mu_{2}\ar[d]^{d^{1}}&\vdots&\vdots&\cdots\\
&& R_{F}[Z_{1}]\ar[d]^{\lambda_{1}}&R_{F}[Z_{1}]\otimes\mu_{2}\ar[d]^{\lambda_{1}\otimes Id}&0&R_{F}[Z_{1}]\otimes \mathbb{Z}/2\mathbb{Z}\ar[d]^{d^{1}}&\cdots\\
&& R_{F}\ar[d]^{\varepsilon}&R_{F}\otimes\mu_{2}\ar[d]^{\varepsilon\otimes H_{1}(\mathrm{inc})}&0&R_{F}\otimes\mathbb{Z}/2\mathbb{Z}\ar[d]^{d^{1}}&\cdots\\
&& \mathbb{Z}\ar[d]^{0}&H_{1}(T,\mathbb{Z})\ar[d]^{C_{\omega}-Id}&H_{2}(T,\mathbb{Z})\ar[d]^{C_{\omega}-Id}&H_{3}(T,\mathbb{Z})\ar[d]^{C_{\omega}-Id}&\cdots\\
&& \mathbb{Z}&H_{1}(T,\mathbb{Z})& H_{2}(T,\mathbb{Z})&H_{3}(T,\mathbb{Z})&\cdots\\
\ar@{-}[rrrrrr]&& &&&&\\
&& &&&&\\}}$
\end{center}

where $C_{\omega}$ denotes the map induced in homology by conjugation $\omega:=\left(
                                                                                 \begin{array}{cc}
                                                                                   0 & -1 \\
                                                                                   1 & 0 \\
                                                                                 \end{array}
                                                                               \right)
$ and $\varepsilon:R_{F}\rightarrow\mathbb{Z}$ is the augmentation map.
\end{lema}

\begin{flushright}
    $\Box$
\end{flushright}

\subsection{Third homology of $SL_{2}$ of a real quadratically closed field}

Since $\lambda_{1}=0$, then $E^{2}_{1,2}=\mathrm{Ker}(\varepsilon)\otimes\mu_{2}=\mathcal{I}_{\mathbf{R}}\otimes\mu_{2}\cong\mu_{2}$ (since $\mathcal{I}_{\mathbf{R}}=\mathbb{Z}\langle\langle-1\rangle\rangle\cong\mathbb{Z}$).

Now $H_{1}(T,\mathbb{Z})=T\cong\mathbf{R}^{\times}$ and the map $C_{\omega}$ induces  $x\mapsto x^{-1}$, so $d^{1}=C_{\omega}-\mathrm{Id}:E^{1}_{1,1}=H_{1}(T,\mathbb{Z})\cong \mathbf{R}^{\times}\rightarrow \mathbf{R}^{\times}\cong H_{1}(T,\mathbb{Z})=E^{1}_{1,0}$ is given by $x\mapsto x^{-2}$. Therefore $E^{2}_{1,0}=\mathbf{R}^{\times}/(\mathbf{R}^{\times})^{2}=\mu_{2}$. Also $\mathrm{Ker}(C_{\omega}-\mathrm{Id})=\mu_{2}=\mathrm{Im}(\varepsilon\otimes H_{1}(\mathrm{inc}))$. It follows that $E^{2}_{1,1}=0$.

Now we have that $C_{\omega}$ induces the identity in $H_{2}(T,\mathbb{Z})=\bigwedge^{2}_{\mathbb{Z}}(\mathbf{R}^{\times})$, thus $d^{1}=C_{\omega}-\mathrm{Id}:E^{1}_{2,1}=H_{2}(T,\mathbb{Z})\cong \mathbf{R}^{\times}\wedge \mathbf{R}^{\times}\rightarrow \mathbf{R}^{\times}\wedge\mathbf{R}^{\times}\cong H_{2}(T,\mathbb{Z})=E^{1}_{2,0}$ is the zero map. Hence $E^{2}_{2,1}=E^{2}_{2,0}=\bigwedge^{2}_{\mathbb{Z}}(\mathbf{R}^{\times})$.

By the arguments of \cite[Section 4.5]{ArticleBlochWignerComplex}, we have that $E^{2}_{3,0}=\mathrm{Tor}^{\mathbb{Z}}_{1}(\mu_{\mathbf{R}},\mu_{\mathbf{R}})$.

Therefore the relevant part of the $E^{2}$-page has the form

\begin{center}
$\scalebox{0.84}{\xymatrix{
&\ar@{-}[ddddd]& \mathcal{RP}_{1}(\mathbf{R})\ar[ddr]^{d^{2}}&E^{2}_{1,3}\ar[ddr]^{d^{2}}&0&\vdots\\
&& \ast&\mu_{2}\ar[ddr]^{d^{2}}&0&\vdots\\
&& \ast&0&\bigwedge^{2}_{\mathbb{Z}}(\mathbf{R}^{\times})&\vdots\\
&& \ast&\ast& \bigwedge^{2}_{\mathbb{Z}}(\mathbf{R}^{\times}) &\mathrm{Tor}^{\mathbb{Z}}_{1}(\mu_{\mathbf{R}},\mu_{\mathbf{R}})\\
\ar@{-}[rrrrr]&& &&&\\
&& &&&\\}}$
\end{center}


\subsection{The calculation of $H_{3}(SL_{2}(\mathbf{R}),\mathbb{Z})$}

Now, $E^{\infty}_{1,2}$ is a subquotient of $E^{2}_{1,2}=\mathcal{I}_{\mathbf{R}}\otimes\mu_{2}\cong\mu_{2}$. So $E^{\infty}_{1,2}=\{1\}$ or $\mu_{2}$.

In \cite[Section 4.7]{ArticleBlochWignerComplex}, it is shown that $E^{\infty}_{0,3}=\mathcal{RB}(F)$  and $E^{\infty}_{0,3}=\mathrm{Tor}^{\mathbb{Z}}_{1}(\mu_{F},\mu_{F})$ for any infinite field $F$. Thus $E^{\infty}_{0,3}=\mathcal{B}(\mathbf{R})$ and $E^{\infty}_{0,3}=\mathrm{Tor}^{\mathbb{Z}}_{1}(\mu_{\mathbf{R}},\mu_{\mathbf{R}})$ respectively.

By Corollary \ref{corollaryrealquadraticallyclosedfields}, we have that $E^{2}_{2,1}=\bigwedge^{2}_{\mathbb{Z}}(\mathbf{R}^{\times})=2\cdot\bigwedge^{2}_{\mathbb{Z}}(\mathbf{R}^{\times})=2\cdot E^{2}_{2,1}$. Now in the section $4.7$ of \cite{ArticleBlochWignerComplex}, it was proved that $2\cdot E^{\infty}_{2,1}=0$. Since $E^{2}_{2,1}$ maps onto $E^{\infty}_{2,1}$, it follows that $E^{\infty}_{2,1}=2\cdot E^{\infty}_{2,1}=0$.

Now let $K:=\mathrm{Ker}(H_{3}(SL_{2}(\mathbf{R}),\mathbb{Z})\rightarrow \mathcal{B}(\mathbf{R}))$. Thus the convergence of the spectral sequence gives us short exact sequences

\begin{center}
$\xymatrix{0\ar[r]&K \ar[r]&H_{3}(SL_{2}(\mathbf{R}),\mathbb{Z})\ar[r]& \mathcal{B}(\mathbf{R})\ar[r]&0     }$
\end{center}

and

\begin{center}
$\xymatrix{0\ar[r]&\mathrm{Tor}^{\mathbb{Z}}_{1}(\mu_{\mathbf{R}},\mu_{\mathbf{R}})  \ar[r]&K\ar[r]& E^{\infty}_{1,2}\ar[r]&0     }$.
\end{center}

In particular, it follows that $|K|=2$ or $4$.


\section{The indecomposable $K_{3}$ of a real quadratically closed field}

\begin{teo}\label{IsommorphismH_3andK_3Ind}Let $\mathbf{R}$ be a real quadratically closed field, then the natural map $H_{3}(SL_{2}(F),\mathbb{Z})\rightarrow K_{3}^{\mathrm{ind}}(F)$ is an isomorphism.
\end{teo}

$\textbf{\emph{Proof}}.$  From the previous remarks and the previous section, we have the following commutative diagram:

\begin{center}
$\xymatrix{0\ar[r]&K\ar[r] \ar@{->>}[d]& H_{3}(SL_{2}(\mathbf{R}),\mathbb{Z}) \ar[r]\ar@{->>}[d]& \mathcal{B}(\mathbf{R})\ar[r]\ar[d]^{=}&0 \\
0\ar[r]&\mathrm{Tor}_{1}^{\mathbb{Z}}(\widetilde{\mu_{\mathbf{R}},\mu_{\mathbf{R}}})=\frac{\mathbb{Z}}{4\mathbb{Z}}\ar[r]&K_{3}^{\mathrm{ind}}(\mathbf{R})\ar[r]&\mathcal{B}(\mathbf{R})\ar[r]&0}$,
\end{center}

The exactness of the bottom row is the main theorem in Suslin (Theorem \ref{Suslinmaintheorem} above). The middle vertical arrow is a surjection by \cite[Lemma 5.1]{HutchinsonLiqunTao}. Thus the map $K\rightarrow\mathrm{Tor}_{1}^{\mathbb{Z}}(\widetilde{\mu_{\mathbf{R}},\mu_{\mathbf{R}}})$ is surjective, and since $|K|\leq4$, it follows that it is an isomorphism and the result follows.

\begin{flushright}
    $\Box$
\end{flushright}

\begin{Rem}
From the proof of Corollary \ref{IsommorphismH_3andK_3Ind}, it follows that the term $E^{\infty}_{1,2}$ of the spectral sequence of above is equal to $\mathcal{I}_{\mathbf{R}}\otimes\mu_{2}\cong\mu_{2}$.
\end{Rem}


Now, let us recall the definition of the Milnor $K$-theory.

For a field $F$, let $T(F^{\times})=\displaystyle\bigoplus_{n\geq0}(F^{\times})^{\otimes n}$, where $(F^{\times})^{\otimes 0}=\mathbb{Z}$, be the tensor algebra over the $\mathbb{Z}$-module $F^{\times}$. Let $J$ be the two-sided homogeneous ideal in $T(F^{\times})$ generated by elements $a\otimes(1-a)$ with $a\in F^{\times}\setminus\{1\}$.

\begin{defi}
\emph{The Milnor K-theory} of a field $F$ is defined to be the graded ring $K_{\bullet}^{M}(F)=T(F^{\times})/J$.
\end{defi}

The image of an element $a_{1}\otimes a_{2}\otimes\cdots\otimes a_{n}$ in $K_{n}^{M}(F)$ is denoted by $\{a_{1},\ldots,a_{n}\}$.

It follows that for any real quadratically closed field $\mathbf{R}$ that $K_{n}^{M}(\mathbf{R})=\langle\{-1,\ldots,-1\rangle\oplus K_{n}^{M}(\mathbf{R})^{0}$ where $K_{n}^{M}(\mathbf{R})^{0}:=\langle\{x_{1},\ldots,x_{n}\}\;| x_{i}>0\; \forall i\rangle$ is $2$-divisible.

Note that, it follows that $K_{n}^{M}(\mathbf{R})^{0}=2\cdot K_{n}^{M}(\mathbf{R})$ when $\mathbf{R}$ is real quadratically closed.


\begin{prop}\label{lemainjectivitythirdhomology}
For any real quadratically closed field $\mathbf{R}$ and any $n\geq3$, the inclusion $SL_{2}(\mathbf{R})\rightarrow SL_{n}(\mathbf{R})$, induces an injection $H_{3}(SL_{2}(\mathbf{R}),\mathbb{Z})\rightarrow H_{3}(SL_{n}(\mathbf{R}),\mathbb{Z})$ and $H_{3}(SL_{n}(\mathbf{R}),\mathbb{Z})\cong H_{3}(SL_{2}(\mathbf{R}),\mathbb{Z})\oplus K_{n}^{M}(\mathbf{R})^{0}$.
\end{prop}

$\textbf{\emph{Proof}}.$ By definition, the map $H_{3}(SL_{2}(\mathbf{R}),\mathbb{Z})\rightarrow K_{3}^{\mathrm{ind}}(\mathbf{R})$ factors through $H_{3}(SL_{3}(\mathbf{R}),\mathbb{Z})$ and hence
$H_{3}(SL_{2}(\mathbf{R}),\mathbb{Z})\rightarrow H_{3}(SL_{3}(\mathbf{R}),\mathbb{Z})$ is injective.

Furthermore $H_{3}(SL_{3}(\mathbf{R}),\mathbb{Z})\cong H_{3}(SL_{n}(\mathbf{R}),\mathbb{Z})$ $\forall n\geq3$ by \cite[Theorem 4.7]{HutchinsonLiqunTao}.

For any infinite field $F$ the sequence

\begin{center}
$\xymatrix{H_{3}(SL_{2}(F),\mathbb{Z})\ar[r]&H_{3}(SL_{3}(F),\mathbb{Z})\ar[r]&2\cdot K_{3}^{M}(F)\ar[r]&0}$
\end{center}

is exact by \cite[Theorem 4.7]{HutchinsonLiqunTao}. Thus, for real quadratically closed fields, we have an exact sequence

\begin{center}
$\xymatrix{0\ar[r]&H_{3}(SL_{2}(\mathbf{R}),\mathbb{Z})\ar[r]&H_{3}(SL_{3}(\mathbf{R}),\mathbb{Z})\ar[r]& K_{3}^{M}(\mathbf{R})^{0}\ar[r]&0}$.
\end{center}

The above sequence is split, since we have a commutative triangle

\begin{center}
$\xymatrix{H_{3}(SL_{2}(\mathbf{R}),\mathbb{Z})\ar[r]\ar[dr]_{\cong}& H_{3}(SL_{3}(\mathbf{R}),\mathbb{Z}) \ar@{->>}[d]\\
&K_{3}^{\mathrm{ind}}(\mathbf{R})}$
\end{center}

where the diagonal arrow is an isomorphism.

\begin{flushright}
    $\Box$
\end{flushright}


\section*{Acknowledgement}
I thank my supervisor Dr. Kevin Hutchinson for comments and advice on this paper. The work leading to this paper was supported by the Research Demonstratorships programme of University College Dublin.


\end{document}